\documentclass[10pt]{article}
\usepackage{amsmath,amssymb,amsfonts}

\newtheorem{theorem}{Theorem}

\date{}

\title{Sampling Part Sizes of Random Integer Partitions}
\bigskip

\author{{\bf Ljuben Mutafchiev}\\
American University in Bulgaria, 2700 Blagoevgrad, Bulgaria \\ and
Institute of Mathematics and Informatics of the \\ Bulgarian
Academy of Sciences
\\ \tt {e-mail: ljuben@aubg.bg; tel.: +359 73 888498}}

\begin{document}
\maketitle

\begin{abstract}
 Let $\lambda$ be a partition of the positive integer $n$, selected
 uniformly at random among all such partitions. Corteel et al.
 (1999) proposed three different procedures of sampling parts of
 $\lambda$ at random. They obtained limiting distributions of the
 multiplicity of the randomly-chosen part as $n\to\infty$. This motivated us to study
 the asymptotic behavior of the part size under the same sampling conditions.
 A limit theorem whenever the part is selected uniformly at random
 among all parts of $\lambda$ (i.e., without any size bias) was
 proved earlier by Fristedt (1993). We consider the remaining two
 (biased) procedures and show that in each of them
 the randomly-chosen part size, appropriately normalized, converges in distribution
 to a continuous random variable. It turns out that different sampling procedures lead to different
 limiting distributions.
\end{abstract}

\vspace {.5cm}

\vspace{.5cm}


\vspace{.2cm}

\section{Introduction}

Partitioning integers into summands (parts) is a subject of
intensive research in combinatorics, number theory and statistical
physics. If $n$ is a positive integer, then by a partition,
$\lambda$, of $n$, we mean a representation
\begin{equation}\label{partition}
 \lambda: \quad n=\sum_{j=1}^n jm_j,
\end{equation}
 in which $m_j$, called multiplicities of parts $j,
j=1,2,...,n$, are non-negative integers. We use $\Lambda(n)$ to
denote the set of all partitions of $n$ and let
$p(n)=\mid\Lambda(n)\mid$. The number $p(n)$ is determined
asymptotically by the famous partition formula of Hardy and
Ramanujan [8]:
 \begin{equation}\label{hardy}
 p(n)\sim\frac{1}{4n\sqrt{3}}
\exp{\left(\pi\sqrt{\frac{2n}{3}}\right)}, \quad n\to\infty.
\end{equation}
A precise asymptotic expansion for $p(n)$ was found later by
Rademacher [13] (more details may be also found in [2; Chapter
5]). Further on, we assume that, for fixed integer $n\ge 1$, a
partition $\lambda\in\Lambda(n)$ is selected uniformly at random
(uar). In other words, we assign the probability $1/p(n)$ to each
$\lambda\in\Lambda(n)$. In this way, each numerical characteristic
of $\lambda$ can be regarded as a random variable defined on
$\Lambda(n)$.

Corteel et al. [3] proposed and studied three procedures of
sampling parts of a random partition $\lambda\in\Lambda(n)$. They
focused on the multiplicities $\mu_{n,j}=\mu_{n,j}(\lambda),
j=1,2,3,$ of the randomly-selected parts and found their limiting
distributions, as $n\to\infty$, in these three cases of sampling.

Our aim in this paper is to determine asymptotically, as
$n\to\infty$, the distributions of the sizes
$\sigma_{n,j}=\sigma_{n,j}(\lambda)$ of the randomly-chosen parts
for the same sampling procedures. One of these three limiting
distributions is obtained earlier by Fristedt [5]. We present the
proofs of the other two limit theorems in Sections 4 and 6. They
combine probabilistic and analytical techniques, which are briefly
described in Section 3. The main results are also stated there.
Sections 2 and 4 contain some definitions, notations and auxiliary
facts.

\section{Preliminaries}
\setcounter{equation}{0}

 We start with the notation $g(x)$ for the generating function of
 the sequence $\{p(n)\}_{n\ge 1}$. For $\mid x\mid<1$, $g(x)$
 admits the well known representation
 \begin{equation}\label{euler}
 g(x)=1+\sum_{n=1}^\infty p(n)x^n =\prod_{k=1}^\infty (1-x^k)^{-1}
 \end{equation}
 (see e.g. [2; Theorem 1.1]).

 For any $\lambda\in\Lambda(n)$ selected uar, we define the random
 variables
 \begin{equation}\label{alph}
 \alpha_j^{(n)} =\alpha_j^{(n)}(\lambda) =\mbox{the number of
 parts of size } j \mbox { in } \lambda
 \end{equation}
and
 \begin{equation}\label{bet}
\beta_j^{(n)}= \beta_j^{(n)}(\lambda) =\left\{\begin{array}{ll} 1 & \qquad  \mbox {if}\qquad \alpha_j^{(n)}>0, \\
 0 & \qquad \mbox {if}\qquad \alpha_j^{(n)}=0.
 \end{array}\right.
 \end{equation}
 Then, obviously, $Z_n=\sum_{j=1}^n \alpha_j^{(n)}$ equals the
 total number of parts and $Y_n=\sum_{j=1}^n \beta_j^{(n)}$ - the
 number of distinct parts in $\lambda\in\Lambda(n)$. Furthermore,
 for $s$ real and $\ge 1$, we also let
 \begin{equation}\label{ex}
 X_{s,n} =\left\{\begin{array}{ll} \sum_{1\le j\le s} j\alpha_j^{(n)} & \qquad  \mbox {if}\qquad 1\le s\le n, \\
 n & \qquad \mbox {if}\qquad s>n;
 \end{array}\right.
 \end{equation}
 \begin{equation}\label{zi}
 Z_{s,n} =\left\{\begin{array}{ll} \sum_{1\le j\le s} \alpha_j^{(n)} & \qquad  \mbox {if}\qquad 1\le s\le n, \\
 Z_n & \qquad \mbox {if}\qquad s>n;
 \end{array}\right.
 \end{equation}
 \begin{equation}\label{uai}
 Y_{s,n} =\left\{\begin{array}{ll} \sum_{1\le j\le s} \beta_j^{(n)} & \qquad  \mbox {if}\qquad 1\le s\le n, \\
 Y_n & \qquad \mbox {if}\qquad s>n.
 \end{array}\right.
 \end{equation}

 Next, by $\mathbb{E}(C)$ we denote the expected value of a random
 variable $C$ with respect to the uniform probability measure $Pr(.)$ defined on the integer partition space $\Lambda(n)$.
 The following two asymptotic equivalences are well known:
 \begin{equation}\label{euai}
 \mathbb{E}(Y_n)\sim\frac{\sqrt{6n}}{\pi}, \quad n\to\infty,
 \end{equation}
 \begin{equation}\label{ezi}
 \mathbb{E}(Z_n)\sim\frac{\sqrt{6n}}{2\pi}\log{n}, \quad n\to\infty
 \end{equation}
 ((\ref{euai}) was proved by Wilf in [15]; a proof of (\ref{ezi}) can
 be found in [3]).

 We introduce now the joint probability generating functions
 \begin{eqnarray}\label{pgf}
 & & \mathbb{E}(y_1^{\alpha_1^{(n)}}y_2^{\alpha_2^{(n)}}
 ...y_n^{\alpha_n^{(n)}})= \nonumber \\
 & & \sum_{\sum jm_j=n, m_j\in\mathbb{N}_0}
 Pr(\alpha_1^{(n)}=m_1,\alpha_2^{(n)}=m_2,...,\alpha_n^{(n)}=m_n)
 \nonumber \\
 & & \times y_1^{m_1}y_2^{m_2}...y_n^{m_n},
 \end{eqnarray}
 where $\mathbb{N}_0=\{0,1,...\}$.  General enumeration methods imply that
 \begin{equation}\label{gf}
 1+\sum_{n=1}^\infty x^n p(n)
 \mathbb{E}(y_1^{\alpha_1^{(n)}}y_2^{\alpha_2^{(n)}}
 ...y_n^{\alpha_n^{(n)}})
 =\prod_{k=1}^\infty \sum_{m_k\in\mathbb{N}_0} (y_k x^k)^{m_k}
 \end{equation}
 (for more details and proof of (\ref{gf}), see [14; Chapter
 V.5]). Here $x$ and $y_k$ are formal
 variables. Setting in (\ref{pgf}) and (\ref{gf}) $y_k=y^k$ for $1\le k\le [s]$
 and $y_k=1$ for $k>[s]$ ($[s]$ denotes
 the integer part of $s$), from
 (\ref{euler}) and (\ref{ex}) we obtain
 $$
 1+\sum_{n=1}^\infty p(n)x^n \mathbb{E}(y^{X_{s,n}})
 =g(x)\prod_{1\le j\le s} \frac{1-x^j}{1-(xy)^j}.
 $$
 A differentiation with respect to $y$ in this identity leads to
 the formula for the expectations we want:
 \begin{equation}\label{expectation}
 \sum_{n=1}^\infty p(n)x^n\mathbb{E}(X_{s,n}) =g(x)h_s(x),
 \end{equation}
 where
 \begin{equation}\label{hs}
 h_s(x)=\sum_{1\le j\le s}\frac{jx^j}{1-x^j}.
 \end{equation}
 Similarly, using (\ref{bet}), (\ref{uai}) and an argument
 developed by Wilf [15], one can deduce that
 $$
 1+\sum_{n=1}^\infty p(n)x^n \mathbb{E}(y^{Y_{s,n}})
 =g(x)\prod_{1\le j\le s}(1+(y-1)x^j)
 $$
 and, in the same way, that
 $$
\sum_{n=1}^\infty p(n)x^n\mathbb{E}(Y_{s,n})
=g(x)\frac{1-x^{[s]+1}}{1-x}.
 $$
It is shown in [11] that, for $s=t\sqrt{6n}/\pi, 0<t<\infty,$
\begin{equation}\label{euaies}
\mathbb{E}(Y_{s,n})\sim (1-e^{-t})\frac{\sqrt{6n}}{\pi}, \quad
n\to\infty.
\end{equation}
Finally, we notice that the asymptotic behavior of $Z_{s,n}$,
defined by (\ref{zi}), is studied in detail by Fristedt [5]. He
assumed that $s/\sqrt{n}\to 0$ as $n\to\infty$ and showed that
$Z_{s,n}$, appropriately normalized, converges in distribution to
a doubly exponential (extreme value) distributed random variable.
A law of large numbers for $Z_{s,n}$ is also proved.

We now proceed with the description of the sampling procedures
introduced by Corteel et al. [3]. We remind that they are three
two-step procedures that combine the outcomes of two experiments.
Therefore, they lead to three different product  probability
spaces (see e.g. [7; Chapter 1.6]). Since in each procedure we
first sample uar a partition $\lambda\in\Lambda(n)$, the
probability space on $\Lambda(n)$, defined in the Introduction, is
included in each product space. The second steps of sampling are,
however, different and therefore, for each different procedure we
obtain a different product space and different product probability
measure. In what follows next we adopt the common notation
$\mathbb{P}(.)$ for the product probability measure of each
sampling procedure, use notation $\sigma_{n,j}$ for the size of
the randomly-chosen part in procedure $j$ ($j=1,2,3$) and follow
the concept of a product space developed in [7; Chapter 1.6]. We
describe the procedures in terms of events
$\{\lambda\in\Lambda(n)\}$, $\{\sigma_{n,j}\le s\}$ ($s\ge 1$) and
their set product $\{\lambda\in\Lambda(n)\}\times\{\sigma_{n,j}\le
s\}$.

{\it Procedure 1.} Given a partition $\lambda\in\Lambda(n)$ chosen
uar (step 1), we select a part uar among all $Y_n$ different parts
of $\lambda$ (step 2). Hence by the product measure formula [7;
Chapter 1.6] and (\ref{uai}),
$$
\mathbb{P}(\{\lambda\in\Lambda(n)\}\times\{\sigma_{n,1}\le s\})
=Pr(\lambda\in\Lambda(n))\left(\frac{Y_{s,n}}{Y_n}\right)
=\left(\frac{1}{p(n)}\right)\left(\frac{Y_{s,n}}{Y_n}\right).
$$
Summation over all $\lambda\in\Lambda(n)$ yields
\begin{equation}\label{procone}
\mathbb{P}(\sigma_{n,1}\le s)
=\mathbb{E}\left(\frac{Y_{s,n}}{Y_n}\right).
\end{equation}

{\it Procedure 2.} Given a partition $\lambda\in\Lambda(n)$ chosen
uar (step 1), we select a part of it with the probability
proportional to its size and multiplicity (step 2). Recalling
definition (\ref{ex}) of the random variable $X_{s,n}$, we obtain
in a similar way that
$$
\mathbb{P}(\{\lambda\in\Lambda(n)\}\times\{\sigma_{n,2} \le s\})
=\left(\frac{1}{p(n)}\right)\left(\frac{X_{s,n}}{n}\right)
$$
and
\begin{equation}\label{proctwo}
\mathbb{P}(\sigma_{n,2}\le s) =\frac{1}{n}\mathbb{E}(X_{s,n}).
\end{equation}

{\it Procedure 3.} Given a partition $\lambda\in\Lambda(n)$ chosen
uar (step 1), we select a part of it uar among all $Z_n$ parts in
$\lambda$ (without any bias, step 2), i.e.,
\begin{equation} \label{procthree}
\mathbb{P}(\{\lambda\in\Lambda(n)\}\times\{\sigma_{n,3}\le s\})
=\left(\frac{1}{p(n)}\right)\left(\frac{Z_{s,n}}{Z_n}\right),
\end{equation}
where $Z_{s,n}$ is defined by (\ref{zi}), and thus
$$
\mathbb{P}(\sigma_{n,3}\le s)
=\mathbb{E}\left(\frac{Z_{s,n}}{Z_n}\right).
$$

{\it Remark.} Each partition $\lambda\in\Lambda(n)$ has a unique
graphical representation called Ferrers diagram [2; Chapter 1.3].
It is obtained as follows. We use the notation $\lambda_k$ to
denote the $k$th largest part of $\lambda$ for $k$ a positive
integer; if the number of parts $Z_n$ of $\lambda$ is $<k$, then
$\lambda_k=0$. The Ferrers diagram illustrates (\ref{partition})
by a two-dimensional array of dots, composed by $\lambda_1$ dots
in the first (most left) row, $\lambda_2$ dots in the second row,
..., $\lambda_{Z_n}$ dots in the last $Z_n$th row. Therefore, a
Ferrers diagram may be considered as a union of disjoint blocks
(rectangles) of dots with base $j$ and height $\alpha_j^{(n)}$
(the multiplicity of part $j$). In this way, the sampling
probability in Procedure 2 is proportional to the area of the
block to which the chosen part belongs.

 In order to make a comparison between the asymptotic behavior of the typical part size $\sigma_{n,j}$ (see next section) and
 its typical multiplicity $\mu_{n,j}$ for procedure $j$ ($j=1,2,3$),
  we also provide the reader with some results, obtained by Corteel et al. [3]:
 \begin{equation}\label{zipf}
 \lim_{n\to\infty}\mathbb{P}(\mu_{n,1}=m) =\frac{1}{m(m+1)}, \quad
 m\ge 1,
 \end{equation}
 \begin{equation}\label{deb}
 \lim_{n\to\infty}\mathbb{P}(\mu_{n,2}=m) =\frac{6(2m+1)}{\pi^2
 m(m+1)^2},\quad m\ge 1
 \end{equation}
 and
 \begin{equation}\label{limit}
\lim_{n\to\infty}\mathbb{P}\left(\frac{2\log{\mu_{n,3}}}{\log{n}}\le
t\right)=t, \quad 0<t<1.
 \end{equation}

\section{Statement of the Main Results and Brief Description of the Method of Proof}
 \setcounter{equation}{0}

Among other important results on random integer partitions
Fristedt [5; p. 712] has proved the following limit theorem.

\begin{theorem}
Let $0<t<1$. Then, we have
$$
\lim_{n\to\infty}\mathbb{P}\left(\frac{2\log{\sigma_{n,3}}}{\log{n}}\le
t\right)=t.
$$
\end{theorem}

{\it Remark.} A comparison between this result and (\ref{limit})
shows that both $\mu_{n,3}$ and $\sigma_{n,3}$, normalized in the
same manner, have one and the same limiting distribution as
$n\to\infty$. Moreover, these results imply that the proportion of
part sizes not greater than $n^{t/2}$ and the proportion of parts
whose multiplicity is $\le n^{t/2}$ are both approximately equal
to $t\in (0,1)$ as $n\to\infty$.

The main results of this paper are devoted to sampling procedures
1 and 2. In the next sections we prove the following limit
theorems for the randomly-chosen part size $\sigma_{n,j}, j=1,2,$
of a random integer partition.

\begin{theorem}
Let $0<t<\infty$. Then, we have
$$
\lim_{n\to\infty}\mathbb{P}\left(\frac{\pi\sigma_{n,1}}{\sqrt{6n}}
\le t\right) =1-e^{-t}.
$$
\end{theorem}

\begin{theorem}
Let $0<t<\infty$. Then, we have
$$
\lim_{n\to\infty}\mathbb{P}\left(\frac{\pi\sigma_{n,2}}{\sqrt{6n}}
\le t\right) =\frac{6}{\pi^2}\int_0^t\frac{u}{e^u-1}du.
$$
\end{theorem}

{\it Remark.} The limiting distribution in Theorem 3 is expressed
in terms of a Debye function (see e.g. [1; Section 27.1]). We
recall that
$$
\int_0^\infty\frac{u}{e^u-1}du =\frac{\pi^2}{6}.
$$
We also notice that Theorems 2 and 3 show that, for sampling
procedures 1 and 2, the typical size of the randomly-chosen part
is of order $const.\sqrt{n}$, while (\ref{zipf}) and (\ref{deb})
imply that its multiplicity is much smaller and approaches a
discrete random variable as $n\to\infty$.

The rest of the paper contains the proofs of Theorems 2 and 3. In
the proof of Theorem 2 we use Fristedt's conditioning device [5],
which allows us to transfer probability distributions of linear
combinations of the multiplicities $\alpha_j^{(n)}$ into
conditional distributions of the corresponding linear combinations
of independent geometrically distributed random variables. This
allows us to approximate the expectations in (\ref{procone}) and
(\ref{procthree}) by the ratios of the expected values of the
corresponding random variables. The proof of Theorem 3 is based on
a Cauchy integral stemming from (\ref{expectation}),
Hardy-Ramanujan's formula (\ref{hardy}) and Hayman's theorem for
estimating coefficients of admissible power series [9] (see also
[4; Chapter VIII.5]).

\section{Proof of Theorem 2}
\setcounter{equation}{0}

We base our proof on (\ref{procone}) and asymptotic equivalences
(\ref{euai}) and (\ref{euaies}). We follow the argument given in
the proof of Theorem 3 of [3]. To replace the expectation in the
right-hand side of (\ref{procone}) by the ratio
$\mathbb{E}(Y_{s,n})/\mathbb{E}(Y_n)$, we need to study how
unlikely is the event
$$
A_n =\left\{\lambda\in\Lambda(n):\mid\frac{\pi Y_n(\lambda)}
{\sqrt{6n}} -1\mid>\epsilon\right\}, \quad \epsilon>0.
$$
Using Fristedt's method [5], Corteel et al. [3] showed that
\begin{equation}\label{probabeien}
Pr(A_n) \le e^{-c\sqrt{n}}, \quad c=c(\epsilon)>0.
\end{equation}

{\it Remark.} Fristedt's approach [5] is based on the identity
\begin{equation}\label{probabgeom}
Pr(\alpha_j^{(n)}=m_j, j=1,...,n) =Pr\left(\gamma_j=m_j,
j=1,...,n\mid\sum_{j\ge 1} j\gamma_j=n\right),
\end{equation}
where $\{\gamma_j\}_{j\ge 1}$ is a sequence of independent
geometrically distributed random variables, whose distribution is
given by
$$
Pr(\gamma_j=k) =(1-q^j)q^{jk}, \quad k=0,1,...
$$
and $\{m_j\}_{j\ge 1}$ are non-negative integers. Eq.
(\ref{probabgeom}) holds for every fixed $q\in (0,1)$. It is
natural to take $q$ so that $Pr(\sum_{j\ge 1}j\gamma_j=n)$ is as
large as possible. Fristedt's almost optimal choice for $q$ is
$q=e^{-\pi/\sqrt{6n}}$. Then, the bound in (\ref{probabeien}) is
easily obtained using this value of $q$.

Next, we represent the probability in (\ref{procone}) in the
following way
\begin{equation}\label{ind}
\mathbb{P}(\sigma_{n,1}\le s)
=\mathbb{E}\left(\frac{Y_{s,n}}{Y_n}I_{A_n^c}\right)
+\mathbb{E}\left(\frac{Y_{s,n}}{Y_n}I_{A_n}\right),
\end{equation}
where $I_{A_n}$ and $I_{A_n^c}$ denote the indicators of events
$A_n$ and $A_n^c$, respectively. Since, for any $\lambda\in
A_n^c$,
$$
\frac{\pi}{\sqrt{6n}(1+\epsilon)}<\frac{1}{Y_n}<\frac{\pi}{\sqrt{6n}(1-\epsilon)}
$$
if $0<\epsilon<1$, the first summand in (\ref{ind}) is estimated
by
\begin{eqnarray}\label{eienc}
& & \mathbb{E}\left(\frac{Y_{s,n}}{Y_n}I_{A_n^c}\right)
=\frac{\pi}{\sqrt{6n}}(1+O(\epsilon))\mathbb{E}(Y_{s,n}I_{A_n^c})
\nonumber \\
& & =\frac{\pi}{\sqrt{6n}}(1+O(\epsilon)) (\mathbb{E}(Y_{s,n})
-\mathbb{E}(Y_{s,n}I_{A_n})).
\end{eqnarray}
Clearly, with probability $1$, we have $Y_{s,n}\le n$. Hence,
using (\ref{probabeien}), we obtain
 $$
 \mathbb{E}(Y_{s,n}I_{A_n}) =O(nPr(A_n)) =O(ne^{-c\sqrt{n}}).
 $$
 Then, for $s=t\sqrt{6n}/\pi$, by (\ref{euaies}) and
 (\ref{eienc}),
 \begin{eqnarray}
 & & \mathbb{E}\left(\frac{Y_{s,n}}{Y_n}I_{A_n^c}\right)
 =\frac{\pi}{\sqrt{6n}}(1+O(\epsilon))\left(\frac{\sqrt{6n}}{\pi}
 (1-e^{-t})\right)(1+o(1)) \nonumber \\
 & & + O(ne^{-c\sqrt{n}}) =1-e^{-t} +O(\epsilon) +o(1)
 +O(ne^{-c\sqrt{n}}). \nonumber
 \end{eqnarray}
 The second term in the right-hand side of (\ref{ind}) is easily
 estimated using (\ref{probabeien}) since it is not greater than
 $Pr(A_n)$. Consequently, (\ref{ind}) becomes
 $$
 \mathbb{P}(\sigma_{n,1}\le s)
 =\mathbb{P}\left(\frac{\pi\sigma_{n,1}}{\sqrt{6n}}\le t\right)
 =1-e^{-t}+O(\epsilon)+o(1).
 $$
 Letting $n\to\infty$ and then $\epsilon\to 0$, we obtain the
 required result.

 \section{Some Remarks on Meinardus Theorem on Weighted Partitions
 and Hayman Admissibility}
 \setcounter{equation}{0}

 This section presents a brief introduction to the
 analytic combinatorics background needed for the proof of Theorem
 3. The starting point in it is eq. (\ref{proctwo}). It
 requires an asymptotic estimate for $\mathbb{E}(X_{s,n})$. By
 (\ref{expectation}), the coefficient $p(n)\mathbb{E}(X_{s,n})$ of
 $x^n$ can be expressed by a Cauchy integral whose integrand is $g(x)h_s(x)$
 (see also
 (\ref{hs})). Its behavior heavily depends on the analytic properties
 of the partition generating function $g(x)$ whose infinite
 product representation (\ref{euler}) shows that its main
 singularity is at $x=1$ (see [2; Chapter 5]). If the integrand of
 the Cauchy integral was only $g(x)$, then a properly chosen contour of
 integration and a proper asymptotic method (the Hardy-Ramanujan
 circle method or the saddle-point method) would yield Hardy-Ramanujan
 asymptotic formula (\ref{hardy}). It has been
 subsequently generalized in various directions most
notably by Meinardus [10] (see also [2; Chapter 6]) who obtained
the asymptotic of the Taylor coefficients of infinite products of
the form
\begin{equation} \label{product}
\prod_{k=1}^\infty (1-x^k)^{-b_k}
\end{equation}
under certain general assumptions on the sequence of non-negative
numbers $\{b_k\}_{k\ge 1}$. Meinardus approach is based on
considering the Dirichlet generating series
\begin{equation} \label{diofzi}
D(z)=\sum_{k=1}^\infty b_k k^{-z}, \quad z=u+iv.
\end{equation}
Below we briefly describe Meinardus assumptions avoiding their
precise statements as well as some extra notations and concepts.
The first assumption ($M_1$) specifies the domain
$\mathcal{H}=\{z: u\ge -C_0\}, 0<C_0<1,$ in the complex plane, in
which $D(z)$ has an analytic continuation. The second one ($M_2$)
is related to the asymptotic behavior of $D(z)$, whenever $\mid
v\mid\to\infty$. A function of the complex variable $z$ which is
bounded by $O(\mid\Im(z)\mid^{C_1}), 0<C_1<\infty$, in certain
domain of the complex plane is called function of finite order.
Meinardus second condition ($M_2$) requires that $D(z)$ is of
finite order in the whole domain $\mathcal{H}$. Finally, the
Meinardus third condition ($M_3$) implies a bound on the ordinary
generating function of the sequence $\{b_k\}_{k\ge 1}$. It can be
stated in a way simpler than the Meinardus original expression by
the inequality
$$
\sum_{k=1}^\infty b_k e^{-k\omega}\sin^2{(\pi ku)} \ge
C_2\omega^{-\epsilon_1}, \quad 0<\frac{\omega}{2\pi}<\mid
u\mid<\frac{1}{2},
$$
for sufficiently small $\omega$ and some constants $C_2,
\epsilon_1>0$ ($C_2=C_2(\epsilon_1)$) (see [6; p. 310]).

It is
known that Euler partition generating function $g(x)$ (which is
obviously of the form (\ref{product})) satisfies the Meinardus
scheme of conditions ($M_1$)-($M_3$) (see e.g. [2; Theorem 6.3]).

In the asymptotic analysis of the Cauchy integral stemming from
(\ref{expectation}) we apply the saddle-point method using Hayman
admissibility theory [9]; see also [4; Chapter VIII.5]. Hayman
studied a wide class of power series satisfying a set of
relatively mild conditions and established general formulas for
the asymptotic order of their coefficients. To present Hayman's
idea and show how it can be applied in the proof of our Theorem 3,
we need to introduce some auxiliary notations.

We consider here a function $G(x)=\sum_{n=1}^\infty G_n x^n$ that
is analytic for $\mid x\mid<\rho, 0<\rho<\infty$. For $0<r<\rho$,
we let
 \begin{equation}\label{ar}
 a(r)=r\frac{G^\prime(r)}{G(r)},
\end{equation}
\begin{equation}\label{br}
b(r)= r\frac{G^\prime(r)}{G(r)}
+r^2\frac{G^{\prime\prime}(r)}{G(r)}
-r^2\left(\frac{G^\prime(r)}{G(r)}\right).
\end{equation}
In the statement of the Hayman's result we use the terminology
given in [4; Chapter VIII.5]. We assume that $G(x)>0$ for $x\in
(R_0,\rho)\subset (0,\rho)$ and satisfies the following three
conditions.

{\it Capture condition.} $\lim_{r\to\rho} a(r)=\infty$ and
$\lim_{r\to\rho} b(r)=\infty$.

{\it Locality condition.} For some function $\delta=\delta(r)$
defined over $(R_0,\rho)$ and satisfying $0<\delta<\pi$, one has
$$
G(re^{i\theta})\sim G(r)e^{i\theta a(r)-\theta^2 b(r)/2}
$$
as $r\to\rho$, uniformly for $\mid\theta\mid\le\delta(r)$.

{\it Decay condition.}
$$
G(re^{i\theta}) =o\left(\frac{G(r)}{\sqrt{b(r)}}\right)
$$
as $r\to\rho$, uniformly for $\delta(r)\le\theta<\pi$.

{\textbf{Hayman Theorem.} Let $G(x)$ be Hayman admissible function
and $r=r_n$ be the unique solution in the interval $(R_0,\rho)$ of
the equation
\begin{equation}\label{eqar}
a(r)=n.
\end{equation}
Then the Taylor coefficients of $G(x)$ satisfy, as $n\to\infty$,
\begin{equation}\label{hayman}
G_n\sim\frac{G(r_n)}{r_n^n\sqrt{2\pi b(r_n)}}
\end{equation}
 with $b(r_n)$ given
by (\ref{br}).

The proof of Theorem 3 (see next section) is divided into two
parts.

A) Proof of Hayman admissibility for $g(x)$.

B) Obtaining an asymptotic estimate for the Cauchy integral
stemming from (\ref{expectation}).

\section{Proof of Theorem 3}
\setcounter{equation}{0}

{\it Part A.}

First we need to show how Hayman's theorem can be applied to find
the asymptotic behavior of the Taylor coefficients of the
partition generating function $g(x)$. Since in (\ref{euler}) we
have $b_k=1, k\ge 1$, the Dirichlet generating series
(\ref{diofzi}) is $D(z)=\zeta(z)$, where $\zeta$ denotes the
Riemann zeta function. We set in (\ref{ar}) and (\ref{br})
$r=r_n=e^{-d_n}, d_n>0$, where $d_n$ is the unique solution of the
equation
\begin{equation}\label{maineq}
a(e^{-d_n})=n.
\end{equation}
((\ref{maineq}) is an obvious modification of (\ref{eqar}).)
Granovsky et al. [6] showed that the first two Meinardus
conditions imply that the unique solution of (\ref{maineq}) has
the following asymptotic expansion:
\begin{equation}\label{dn}
d_n =\sqrt{\zeta(2)/n} +\frac{\zeta(0)}{2n} +O(n^{-1-\beta})
=\frac{\pi}{\sqrt{6n}}-\frac{1}{4n} +O(n^{-1-\beta}),
\end{equation}
where $\beta>0$ is fixed constant (here we have also used that
$\zeta(0)=-1/2$; see [1; Chapter 23.2]). We also notice that
(\ref{br}) and (\ref{dn}) impliy that
\begin{equation}\label{bdn}
b(e^{-d_n}) =2\zeta(2)d_n^{-3} +O(d_n^{-2}) \sim\frac{\pi^2}{3}
d_n^{-3} \sim\frac{2\sqrt{6}}{\pi} n^{3/2}
\end{equation}
(see [12; Lemma 2.2] with $D(z)=\zeta(z)$). Hence, by
(\ref{maineq}) and (\ref{bdn}), $a(e^{-d_n})\to\infty$ and
$b(e^{-d_n})\to\infty$ as $n\to\infty$, that is, Hayman's
``capture'' condition is satisfied with $r=r_n=e^{-d_n}$. To show
next that Hayman's ``decay'' condition is satisfied by $g(x)$ we
set
\begin{equation}\label{delta}
\delta_n =\frac{d_n^{4/3}}{\Omega(n)}
=\frac{\pi^{4/3}}{(6n)^{2/3}\Omega(n)}
\left(1+O\left(\frac{1}{\sqrt{n}}\right)\right)
\end{equation}
with $d_n$ given by (\ref{dn}), where $\Omega(n)\to\infty$ as
$n\to\infty$ arbitrarily slowly. We can apply now an estimate for
$\mid g(e^{-d_n+i\theta})\mid$ established in a general form in
[12; Lemma 2.4] using all three Meinardus conditions. It states
that there are two positive constants $c_0$ and $\epsilon_0$, such
that, for sufficiently large $n$,
\begin{equation}\label{decay}
\mid g(e^{-d_n+i\theta})\mid \le
g(e^{-d_n})e^{-c_0d_n^{-\epsilon_0}}
\end{equation}
uniformly for $\delta_n\le\mid\theta\mid<\pi$. This, in
combination with (\ref{bdn}), implies that $\mid
g(e^{-d_n+i\theta})\mid=o(g(e^{-d_n})/\sqrt{b(e^{-d_n})})$
uniformly in the same range for $\theta$, which is just Hayman's
``decay'' condition. Finally, by Lemma 2.3 of [12], established
using Meinardus conditions ($M_1$) and ($M_2$), Hayman's
``locality'' condition is also satisfied by $g(x)$. In fact, this
lemma implies in the particular case $D(z)=\zeta(z)$ that
\begin{equation}\label{locality}
e^{-i\theta n}\frac{g(e^{-d_n+i\theta})}{g(e^{-d_n})}
=e^{-\theta^2 b(e^{-d_n})/2}(1+O(1/\Omega^3(n))
\end{equation}
uniformly for $\mid\theta\mid\le\delta_n$, where $b(e^{-d_n})$ and
$\delta_n$ are determined by (\ref{bdn}) and (\ref{delta}),
respectively.
 Hence all conditions of
Hayman's theorem hold and we can apply it with $G_n=p(n),
G(x)=g(x), r_n=e^{-d_n}$ and $\rho=1$ to find that
\begin{equation}\label{asypn}
p(n)\sim\frac{e^{nd_n}g(e^{-d_n})}{\sqrt{2\pi b(e^{-d_n})}}, \quad
n\to\infty.
\end{equation}

{\it Remark.} To show that formula (\ref{asypn}) yields
(\ref{hardy}), one has to replace (\ref{dn}) and (\ref{bdn}) in
the right hand side of (\ref{asypn}).
 The asymptotic of $g(e^{-d_n})$ is
determined by a general lemma due to Meinardus [10] (see also [2;
Lemma 6.1]). Since $\zeta(0)=-1/2$ and
$\zeta^\prime(0)=-\frac{1}{2}\log{(2\pi)}$ (see [1; Chapter
23.2]), in the particular case of $g(e^{-d_n})$ this lemma implies
that
\begin{eqnarray}
& & g(e^{-d_n}) =\exp{(\zeta(2)d_n^{-1} -\zeta(0)\log{d_n}
+\zeta^\prime(0) +O(d_n^{c_1}))} \nonumber \\
& & =\exp{\left(\frac{\pi^2}{6d_n} +\frac{1}{2}\log{d_n}
-\frac{1}{2}\log{(2\pi)} +O(d_n^{c_1})\right)}, \quad n\to\infty,
\nonumber
\end{eqnarray}
 where $0<c_1<1$. The rest of the computation leading to (\ref{hardy}) is based on simple
 algebraic manipulations and cancellations.

 {\it Part B.}

 We are now ready to apply
 Cauchy coefficient formula to (\ref{expectation}). We use the
 circle $x=e^{-d_n+i\theta}, -\pi<\theta\le\pi$, as a contour of
 integration and obtain
  $$
 p(n)\mathbb{E}(X_{s,n}) =\frac{e^{nd_n}}{2\pi}\int_{-\pi}^\pi
 g(e^{-d_n+i\theta}) h_s(e^{-d_n+i\theta}) e^{-i\theta n}d\theta.
 $$
 Then, we break up the range of integration as follows:
 \begin{equation} \label{sumint}
p(n)\mathbb{E}(X_{s,n}) =J_1(s,n) +J_2(s,n),
 \end{equation}
 where
 \begin{equation}\label{jone}
 J_1(s,n) =\frac{e^{nd_n}}{2\pi}\int_{-\delta_n}^{\delta_n}
g(e^{-d_n+i\theta}) h_s(e^{-d_n+i\theta}) d\theta,
\end{equation}
\begin{equation} \label{jtwo}
J_2(s,n)
=\frac{e^{nd_n}}{2\pi}\int_{\delta_n<\mid\theta\mid\le\pi}
g(e^{-d_n+i\theta}) h_s(e^{-d_n+i\theta}) d\theta
\end{equation}
and $\delta_n$ is defined by (\ref{delta}).

To estimate $J_2(s,n)$, for $s=t\sqrt{6n}/\pi, 0<t<\infty$, we
notice that by the definition of Riemann integrals and (\ref{dn}),
\begin{eqnarray}\label{hsest}
& & \mid h_s(e^{-d_n+i\theta})\mid \le \sum_{1\le j\le s}
\frac{je^{-jd_n}}{\mid 1-e^{-jd_n+ij\theta}\mid} \nonumber \\
 & & \le d_n^{-2}\sum_{d_n\le jd_n\le
 sd_n}\frac{jd_ne^{-jd_n}}{1-e^{-jd_n}}d_n \nonumber \\
 & & \sim d_n^{-2}\int_0^t\frac{ue^{-u}}{1-e^{-u}}du=O(d_n^{-2}) =O(n).
\end{eqnarray}
Combining (\ref{dn}), (\ref{bdn}), (\ref{decay}), (\ref{asypn}),
(\ref{jtwo}) and (\ref{hsest}), we obtain
 \begin{eqnarray}\label{jtwoest}
 & & \mid J_2(s,n)\mid \le\frac{e^{nd_n}}{2\pi}\int_{\delta_n<\mid\theta\mid\le\pi}
\mid g(e^{-d_n+i\theta}) h_s(e^{-d_n+i\theta})\mid d\theta
\nonumber \\
 & & =O(e^{nd_n} g(e^{-d_n})ne^{-c_0d_n^{-\epsilon_0}})
 =O\left(\frac{e^{nd_n}g(e^{-d_n})}{\sqrt{b(e^{-d_n})}} n^{1+3/4}
e^{-c_0d_n^{-\epsilon_0}}\right) \nonumber \\
& & =O(p(n)n^{7/4}e^{-c_0d_n^{-\epsilon_0}})
=O(p(n)n^{7/4}e^{-c_2n^{\epsilon_0/2}}) =o(np(n)),
 \end{eqnarray}
where $c_2>0$.

The estimate of $J_1(s,n)$ follows from Hayman's ``locality''
condition (\ref{locality}). We also need to expand $h_s$ by Taylor
formula in the following way:
\begin{eqnarray}\label{hstaylor}
& & h_s(e^{-d_n+i\theta}) =h_s(e^{-d_n})
+O\left(\mid\theta\mid\frac{d}{dx} h_s(x)\mid_{x=e^{-d_n}}\right)
\nonumber \\
& & =h_s(e^{-d_n}) +O\left(\delta_n \frac{d}{dx}
h_s(x)\mid_{x=e^{-d_n}}\right).
\end{eqnarray}
For $s=t\sqrt{6n}/\pi$, we can consider, as previously, the sum
representing $h_s(e^{-d_n})$ as a Riemann sum. So, we can replace
it by the corresponding integral. Thus, by (\ref{hs}) and
(\ref{dn}), we have
\begin{eqnarray}
& & h_s(e^{-d_n}) =d_n^{-2}\sum_{d_n\le jd_n\le t\sqrt{6n}d_n/\pi}
\frac{jd_n e^{-jd_n}}{1-e^{-jd_n}}d_n \nonumber \\
& & \sim d_n^{-2}\int_{d_n}^{\frac{\sqrt{6n}}{\pi}d_n t}
\frac{ue^{-u}}{1-e^{-u}}du =\frac{6n}{\pi^2} \int_0^t
\frac{u}{e^u-1}du +O(\sqrt{n}) \nonumber
 \end{eqnarray}
since $d_n^{-2}=6n/\pi^2+O(\sqrt{n})$. In the same way we can
estimate the first derivative of $h_s$:
\begin{eqnarray}
& & \frac{d}{dx} h_s(x)\mid_{x=e^{-d_n}} =\sum_{1\le j\le s}
\frac{j^2 e^{(j-1)d_n}}{(1-e^{-jd_n})^2}  \sim
d_n^{-3}\sum_{d_n\le jd_n\le t\sqrt{6n}d_n/\pi}
\frac{(jd_n)^2 e^{-jd_n}}{(1-e^{-jd_n})^2}d_n \nonumber \\
& & \sim d_n^{-3}\int_0^t \frac{u^2 e^{-u}}{(1-e^{-u})^2}du
=O(n^{3/2}). \nonumber
\end{eqnarray}
Hence, by (\ref{delta}), the error term in (\ref{hstaylor})
becomes
$$
O\left(\delta_n\frac{d}{dx} h_s(x)\mid_{x=e^{-d_n}}\right)
=O(n^{5/6}/\Omega(n))
$$
and therefore, uniformly for $\mid\theta\mid\le\delta_n$,
$$
h_s(e^{-d_n+i\theta}) =\frac{6n}{\pi^2}\int_0^t\frac{u}{e^u-1}du
+o(n).
$$
Inserting this estimate and (\ref{locality}) into (\ref{jone}) and
applying the asymptotic of the partition function $p(n)$ from
(\ref{asypn}), we obtain
\begin{eqnarray}\label{joneest}
& & J_1(s,n) =\frac{e^{nd_n}g(e^{-d_n})}{2\pi}
\left(\int_{-\delta_n}^{\delta_n} e^{-\theta^2 b(e^{-d_n})/2}
(1+O(1/\Omega^3(n))d\theta\right) \nonumber \\
& & \times\left(\frac{6n}{\pi^2}\int_0^t\frac{u}{e^u-1}du
+o(n)\right) \nonumber \\
& & \sim\frac{e^{nd_n}g(e^{-d_n})}{\sqrt{b(e^{-d_n})}2\pi}
\left(\int_{-\delta_n\sqrt{b(e^{-d_n})}}^{\delta_n\sqrt{b(e^{-d_n})}}
e^{-y^2/2}dy\right)
\left(\frac{6n}{\pi^2}\int_0^t\frac{u}{e^u-1}du\right) \nonumber
\\
& & \sim\frac{e^{nd_n}g(e^{-d_n})}{\sqrt{b(e^{-d_n})}2\pi}
\left(\int_{-\infty}^\infty e^{-y^2/2}dy\right)
\left(\frac{6n}{\pi^2}\int_0^t\frac{u}{e^u-1}du\right) \nonumber
\\
& & =\frac{e^{nd_n}g(e^{-d_n})}{\sqrt{2\pi
b(e^{-d_n})}}\left(\frac{6n}{\pi^2}\int_0^t\frac{u}{e^u-1}du\right)
\nonumber \\
& & \sim
p(n)\left(\frac{6n}{\pi^2}\int_0^t\frac{u}{e^u-1}du\right),
\end{eqnarray}
where for the second asymptotic equivalence we have used
(\ref{bdn}) and (\ref{delta}) in order to get
$$
\delta_n\sqrt{b(e^{-d_n})}\sim
\frac{\pi^{5/6}\sqrt{2}}{6^{1/6}\Omega(n)} n^{1/12}\to\infty
$$
if $\Omega(n)\to\infty$ as $n\to\infty$ not too fast, so that
$\frac{n^{1/12}}{\Omega(n)}\to\infty$. It is now clear that, for
$s=t\sqrt{6n}/\pi$, (\ref{sumint})-(\ref{jtwo}), (\ref{jtwoest})
and (\ref{joneest}) yield
$$
p(n)\mathbb{E}(X_{s,n}) =p(n)\frac{6n}{\pi^2}
\int_0^t\frac{u}{e^u-1}du +o(np(n))
$$
and therefore
$$
\mathbb{E}(X_{s,n})\sim \frac{6n}{\pi^2}
\int_0^t\frac{u}{e^u-1}du.
$$
The result of Theorem 3 follows immediately from (\ref{proctwo}).

\section*{Acknowledgements}

I am grateful the referee for carefully reading the paper and for
his helpful comments.

\end{document}